\documentstyle[amscd,amssymb,verbatim,12pt]{amsart}
\theoremstyle{plain}
\newtheorem{Thm}{Theorem}[section]

\newtheorem{Prop}[Thm]{Proposition}
\newtheorem{Cor}[Thm]{Corollary}

\theoremstyle{definition}
\newtheorem{Rem}[Thm]{Remark}

\newtheorem{Def}[Thm]{Definition}
\numberwithin{equation}{section}
\newcommand{\eqnum}{\setcounter{equation}{\value{Thm}}}
\newcommand{\thnum}{\setcounter{Thm}{\value{equation}}}

\begin{document}
\title[ON NORMAL SUBGROUPS OF PRODUCT OF GROUPS]%
	{ON NORMAL SUBGROUPS OF PRODUCT OF GROUPS}
\author[A. K. DAS]%
	{ ASHISH KUMAR DAS }
\address{Department of Mathematics, North Eastern Hill
University, Permanent Campus, Shillong-793022, Meghalaya, India.} 
\email{akdas@@nehu.ac.in}
\thanks{}
\date{}
\begin{abstract}
The object of this paper is to find a necessary and sufficient condition for
the groups $G_1, G_2, \dots ,G_n$ so that  every normal subgroup of the
product $\prod_{i=1}^{n} G_i$ is of the  type $\prod_{i=1}^{n} N_i$ with $N_i
\trianglelefteq G_i$, $i=1,2, \dots ,n$. As a consequence we obtain a
well-known result due to R. Remak about centreless completely reducible groups
having finitely many direct factors.  \end{abstract} 
\maketitle

\section{Introduction} \label{S:intro}

Let $G$ be a finite group. Then $G$ has a composition series given by $\{e\} =
H_0 \vartriangleleft H_1  \vartriangleleft \dots  \vartriangleleft H_{n-1} 
\vartriangleleft H_n = G.$ The set-with-multiplicities of composition factors
associated to this composition series is given by $C(G) = \{H_i / H_{i-1}
: i = 1,2, \dots , n\};$ noting that any two composition series of $G$ give
rise to the same set-with-multiplicities of composition factors, upto
isomorphism of factors (see \cite{nJ74}, Jordan-Holder Theorem). It is a
standard fact (see \cite{jF82}) that  $\forall \; K  \vartriangleleft  G$, 
$C(G)$ is the disjoint union (i.e. union counting multiplicities) of  
$C(G/K)$ and $C(K)$. By convention $C(\{e\}) = \emptyset$.

 In \cite{tL01}, Leinster has proved that if $G_1$and $G_2$ are any two
finite groups such that  $C(G_1)$ and $C(G_2)$ have no member in common
then every normal subgroup of $G_1 \times G_2$ is of the type $N_1 \times
N_2$ where $N_i \trianglelefteq G_i$, $i = 1,2$. However, considering $G_1 =
G_2 =$ the alternating group $A_5$ (a simple non-abelian group), one can see
that the converse is not true.

In this paper we derive a condition which is necessary as well as sufficient
for any pair of groups $G_1$ and $G_2$ (finite or infinite, abelian or
non-abelian) so that all normal subgroups of $G_1 \times G_2$ are of the above
mentioned type. We then extend our result to any finite number of groups
and obtain, as a consequence, a well-known result of R. Remak about
centreless completely reducible groups having finitely many direct factors.
The motivation for this work lies in finding  examples of  non-abelian perfect
groups;  a perfect group being one in  which the sum of the orders of all
normal subgroups equals twice the order of the group (see \cite{tL01}).

\section{The NS-condition} \label{S:nscon}

We begin this section with a few definitions.

\begin{Def}\label{D:std}
Let $G_1$ and $G_2$ be any two groups. A normal subgroup $N$ of $G_1 \times
G_2$  is said to be of standard type if $N = N_1 \times N_2$ where $N_i
\trianglelefteq G_i$,  $i = 1,2.$
\end{Def}

\begin{Def}\label{D:com}
Two groups $G_1$ and $G_2$ are said to have a subgroup in common if there exist
non-trivial subgroups $H_i$ of $G_i,$ $i= 1,2,$ such that
$H_1 \cong H_2$.
\end{Def}

\begin{Def}\label{D:nsc}
Two groups $G_1$ and $G_2$ are said to satisfy NS-condition if $\forall \; H_i
\vartriangleleft G_i$, $i= 1,2,$ the centres $Z(G_1 /H_1)$ and  $Z(G_2 /H_2)$
of the quotient groups $G_1 /H_1$ and  $G_2 /H_2$ have no subgroup in common. 
\end{Def}

\begin{Rem}\label{R:csg}
If $H \trianglelefteq G$ then the subgroups of $Z(G /H)$ are of the type $K/H$
where $H \trianglelefteq K \leq G $ and  $gkg^{-1}k^{-1} \in H$, \;
$\forall g \in G$ and $\forall k \in K$. In particular  $K \trianglelefteq G.$
\end{Rem}

The following proposition shows that NS-condition is strictly weaker than
Leinster's hypothesis (see Proposition 3.1 of \cite{tL01}).

\begin{Prop}\label{P:com}
If $G_1$ and $G_2$ are two finite groups such that  $C(G_1)$ and
$C(G_2)$ have no abelian member in common then $G_1$ and $G_2$ satisfy
NS-condition. However, the converse is not true.
\end{Prop}

\begin{pf}
Suppose $H_i \vartriangleleft G_i$, $i = 1,2,$ are such that the centres
$Z(G_1 /H_1)$  and  $Z(G_2 /H_2)$ have a subgroup in common. So, there are
non-trivial subgroups $K_i /H_i$ of  $Z(G_i /H_i)$, $i = 1,2,$   such that
$K_1 /H_1 \cong K_2 /H_2$. Then $C(K_1 /H_1) = C(K_2 /H_2)$. Since  $H_i
\vartriangleleft K_i \trianglelefteq G_i,$  we have $C(K_i /H_i)$ $ \subseteq
C(K_i) \subseteq  C(G_i)$, $i=1,2$. Thus,  $C(G_1)$ and $C(G_2)$   have an
abelian member in common; noting that  both $K_1 /H_1$ and $K_2 /H_2$ are
abelian.

That the converse is not true can be easily seen by considering $G_1 = S_4$
and $ G_2 =  {\Bbb{Z}}/ 3{\Bbb{Z}}$, and noting that  $S_4$ and $S_4/V \cong
S_3$ have trivial centres; $V$ being the normal subgroup $\{ e, (12)(34),
(13)(24), (14)(23) \}$ of $S_4$.
\end{pf}

\begin{Prop}
Two finite groups $G_1$ and $G_2$ satisfy NS-condition if and only if
$gcd(n_1 ,n_2) = 1$ where

\[
\textstyle
n_i= \underset{H_i \trianglelefteq G_i}{\prod}
|Z(G_i /H_i)|, \; \; i= 1,2.
\]

\end{Prop}

\begin{pf}

If each $n_i$ is divisible by a prime $p$, then p divides $|Z(G_i /H_i)|$ for
some $H_i \trianglelefteq G_i$,  $i = 1,2$. So, each of $Z(G_i /H_i)$,  $i
= 1,2$,\; has a subgroup of order $p$ (see \cite{jF82}) . Obviously such a
subgroup is isomorphic to ${\Bbb{Z}}/ p{\Bbb{Z}}$.

Conversely, if $gcd(n_1 ,n_2) =1$ then  $\forall \; H_i 
\trianglelefteq G_i$,  $i = 1,2$, \; the centres
$Z(G_1 /H_1)$ and  $Z(G_2 /H_2)$  have coprime orders
and so they can not have any subgroup in common.
\end{pf}

The following corollary is immediate.

\begin{Cor}
If two finite groups  $G_1$ and $G_2$ have coprime orders then they satisfy
NS-condion. Converse also holds if  $G_1$ and $G_2$ are finite abelian.
\end{Cor}

\section{Main Theorem}\label{S:proof}

In this section we show that the NS-condition is indeed necessary as well as
sufficient.

\begin{Thm}
Let $G_1$ and $G_2$ be any two groups (finite or infinite, abelian or
non-abelian). Then all normal subgroups of $G_1 \times G_2$ are of standard
type if and only if $G_1$ and $G_2$ satisfy NS-condition.
\end{Thm}

\begin{pf}
Suppose $G_1$ and $G_2$ satisfy NS-condition. Let $N \trianglelefteq G_1
\times G_2$. Set $H_1 = \pi_1 ((G_1 \times \{ e_2\}) \cap N)$ and $H_2 = \pi_2
((\{e_1\} \times  G_2) \cap N)$ where $e_i$ are identities of $G_i$ and
$\pi_i : G_1 \times G_2 \longrightarrow G_i$ are projections, $i = 1,2$. Then
$H_1 \times \{e_2\}, \{e_1\} \times H_2  \subset N \subset \pi_1 (N) \times
\pi_2 (N)$ and so 
\eqnum
\begin{equation} \label{E:prod}
H_1 \times H_2  = (H_1 \times \{e_2\})( \{e_1\} \times H_2) \subset  \pi_1 (N)
\times \pi_2 (N) \end{equation}
\thnum

\noindent It may be noted here that $H_i \trianglelefteq G_i$ and $H_i
\trianglelefteq \pi_i (N)$, $i=1,2.$  Now, suppose $a_1 \in \pi_1 (N)$. Then
$(a_1, a_2) \in N$ for some $a_2 \in G_2$; in fact $a_2 \in \pi_2 (N)$.
Therefore, $\forall g_1 \in G_1,$ we have 
\begin{align*}
& (g_1 a_1 {g_1}^{-1}, a_2) = (g_1,e_2)(a_1, a_2)({g_1}^{-1},e_2) \in N \\
\Longrightarrow \; & (g_1 a_1 {g_1}^{-1} {a_1}^{-1}, e_2) \in N \\ 
\Longrightarrow \; & g_1 a_1 {g_1}^{-1} {a_1}^{-1} \in H_1  \\
\Longrightarrow \; & g_1 H_1 a_1H_1 = a_1 H_1 g_1 H_1\in G_1.
\end{align*}

Thus, $a_1H_1 \in Z(G_1 /H_1)$.  So, we have $ \pi_1 (N)/H_1 \subset   Z(G_1
/H_1).$  Similarly, $ \pi_2 (N)/H_2 \subset   Z(G_2 /H_2).$  Note that if
$a_1, b_1 \in \pi_1 (N)$ then $(a_1, a_2),(b_1, b_2) \in N$ for some $a_2, b_2
\in \pi_2 (N),$  and so $(a_1 {b_1}^{-1}, a_2 {b_2}^{-1}),$ $ (a_1 b_1, a_2
b_2) \in N$. Therefore,  \begin{align*}
&a_1 H_1 =b_1 H_1 \; \Longleftrightarrow \; a_1 {b_1}^{-1} \in H_1
\; \Longleftrightarrow \; (a_1 {b_1}^{-1}, e_2) \in N \\
\Longleftrightarrow \; & (e_1, a_2 {b_2}^{-1}) \in N \; \Longleftrightarrow \;
a_2 {b_2}^{-1} \in H_2 \; \Longleftrightarrow \; a_2 H_2 =b_2 H_2.
\end{align*}

\noindent This means that we have a well-defined injective map $ f : \pi_1
(N) / H_1 $ $\longrightarrow \pi_2 (N) / H_2$ given by $f(a_1 H_1) = a_2 H_2$
where $(a_1, a_2) \in N.$  Also, $f(a_1 H_1 b_1 H_1) = f(a_1
b_1 H_1) = a_2 b_2 H_2 = a_2 H_2 b_2 H_2 = f(a_1 H_1) f(b_1 H_1)$, showing that
$f$ is a homomorphism.  Finally, if $b \in \pi_2 (N)$ then $(a, b) \in N$ for
some $a \in \pi_1 (N)$ and so $f(a H_1) = b H_2$, which implies that $f$ is
surjective.  Thus $f$ is an isomorphism.  Hence it follows from the hypothesis
that $\pi_i (N) / H_i$  are trivial subgroups of $Z(G_i /H_i)$, $i=1,2$. 
Therefore, $H_i = \pi_i (N),$ $i=1,2$, and so $N= H_1 \times H_2,$ by
(\ref{E:prod}).

Conversely, suppose $G_1$ and $G_2$ do not satisfy the NS-condition. So, there
exist $H_i \vartriangleleft G_i$, $i =1,2$, such that $Z(G_1 /H_1)$ and
$Z(G_2 /H_2)$ have a subgroup in common.  Let $K_i /H_i$ be non-trivial
subgroups of  $Z(G_i /H_i)$, $i=1,2$, such that there is an isomorphism $F :
K_1 /H_1 \longrightarrow  K_2 /H_2.$  Put $N = \{(a_1, a_2) \in K_1 \times K_2
: F(a_1 H_1) = a_2 H_2 \}.$  Let $(a_1, a_2), (b_1, b_2) \in N$ then $F(a_1
H_1) = a_2 H_2$ and $F(b_1 H_1) = b_2 H_2 $.  So, $F(a_1 {b_1}^{-1} H_1) =
a_2 {b_2}^{-1} H_2$.  Thus $(a_1, a_2)(b_1, b_2)^{-1} = (a_1 {b_1}^{-1},
a_2 {b_2}^{-1}) \in N$, showing that $N$ is a subgroup of $G_1 \times G_2$. 
Again let $(a_1, a_2) \in N$ and $(g_1, g_2) \in G_1 \times G_2$. Then, 
$(g_1, g_2) (a_1, a_2) (g_1, g_2)^{-1} = (g_1 a_1 {g_1}^{-1}, g_2 a_2
{g_2}^{-1}) \in K_1 \times K_2$, since, by Remark \ref{R:csg}, $K_i
\trianglelefteq G_i$, $i=1,2$. Also, since $a_i H_i \in K_i /H_i \subset
Z(G_i /H_i)$, $i= 1,2$, we have

\[
F(g_1 a_1 {g_1}^{-1} H_1) = F(a_1 H_1)= a_2 H_2 = g_2 a_2 {g_2}^{-1}
H_2.
\]

\noindent Thus $(g_1, g_2) (a_1, a_2) (g_1, g_2)^{-1} \in N$, and so $N
\trianglelefteq G_1 \times G_2$. On the other hand, suppose $N$ is of standard
form $N_1 \times N_2$ where $N_i \trianglelefteq G_i$, $i=1,2$. Then, $ \pi_i
(N) = N_i$, $i=1,2$.  But since F is bijective, we have $\pi_i (N) =K_i$,
$i=1,2$. Therefore, $N = K_1 \times K_2$. Since $K_1 /H_1$ is non-trivial,
there is some $a_1 \in K_1$ such that $a_1 H_1 \neq H_1$.  But $(a_1, e_2) \in 
K_1 \times K_2 = N$.  So, $F(a_1 H_1) = e_2 H_2 = H_2$, the zero element of
$K_2 /H_2$.  Therefore, since $F$ is injective, we have $a_1 H_1 = H_1$, the
zero element of $K_1 /H_1$.  This contradiction shows that $N$ is not of
standard type.
\end{pf}
 
Note that $(G_1 \times G_2)/(H_1 \times H_2) \cong  G_1/H_1 \times G_2/H_2$
where  $H_i \trianglelefteq G_i$, $ i = 1,2$. Therefore, using induction and
the fact  that the centre of a product of groups equals the product of the
centres of the constituent groups, we have the following corollary:

\begin{Cor}
Let $G_1, G_2, \dots ,G_n$ be any $n$ groups (finite or infinite, abelian or
non-abelian). Then every normal subgroup of the product $\prod_{i=1}^{n}
G_i$ is of the  type $\prod_{i=1}^{n} N_i$ with $N_i \trianglelefteq
G_i$, $i=1,2, \dots ,n$,  if and only if  $G_1, G_2, \dots ,G_n$ satisfy
NS-condition pairwise. 
\end{Cor}

\begin{Rem}
If atleast $n-1$ of the groups  $G_1, G_2, \dots ,G_n$ are
simple non-abelian then they satisfy NS-condition pairwise. In view of this, 
the above corollary says, in particular, that if G is a direct product of
non-abelian simple groups $G_i$, $i=1,2, \dots ,n$, then every normal subgroup
of $G$ is a direct product of certain of the $G_i$'s. This is  a well-known
result proved by R. Remak (see \cite{dR96}, page 88).
\end{Rem}

\end{document}